\documentclass[10pt]{amsart}

\usepackage{amsmath, amsthm, amssymb, amsfonts, amscd}

\usepackage{graphicx}
\usepackage[bookmarks=false]{hyperref}

\newcommand{\HausB}{\dim_{\mathrm{HB}}}
\newtheorem{thm}{Theorem}[section]
\newtheorem{defn}[thm]{Definition}
\newtheorem{rem}[thm]{Remark}
\newtheorem{exm}[thm]{Example}

\date{}

\begin{document}

\title{Fractal Dimension of the Kronecker Product Based Fractals}

\author{Anatoly E. Voevudko}
\address{CSULB - California State University, Long Beach, CA  ... February 22, 2018}
\email{AnatolV@outlook.com}


\begin{abstract}
A simple method of calculating the Hausdorff-Besicovitch dimension of the Kronecker Product based 
fractals is presented together with a compact R script realizing it. 
The proposed new formula is based on traditionally used values of the number of self-similar objects and 
the scale factor that are now calculated using appropriate values of both the initial fractal matrix and the second order resultant matrix. 
This method is reliable and producing dimensions equal to many already determined values of well-known fractals.
\end{abstract}
\maketitle

\section{Introduction}

The Kronecker product [1] generates fractals if it applies to the same matrix [1-3], i.e., it is the \emph{self Kronecker product}.

The origin, nature and generating technique related to the Kronecker product based (KPB) fractals  are explained
in more details in [2,3]. A lot of well-known and new fractal samples are shown and KPB fractals generators
are offered both in JavaScript and in R [2,3].

A simple method of calculating the Hausdorff-Besicovitch dimension (HBD) of the Kronecker Product based 
fractals was discovered and is presented together with a compact R script realizing it. 
The proposed new formula is based on traditionally used values of the number of self-similar objects and 
the scale factor that are now calculated using  appropriate values of both the initial fractal matrix and the second order resultant matrix. 
This method is reliable and producing dimensions equal to many already determined values of well-known canonical fractals.

It should be stressed in the beginning that this method only works and could be applied to any fractal matrix (defined below) 
to determine KPB fractal dimension.

\section{Generating and plotting fractals}

Let's start with a few definitions that are related and applied only to matrices used for fractal images generation using
a \emph{self Kronecker product} (a.k.a. the \emph{Kronecker power}).

\begin{defn}
		\textbf{The Kronecker power} $n$ of a matrix $M$ is defined as\\
		$M \otimes M \otimes ...  \otimes M$, i.e., n times self Kronecker product. In short: $M^{n\otimes}$.
\end{defn}

\begin{rem}
	It is often said that matrix $M$ has an order (or level) $n$ if it has the Kronecker power $n$. Same order is used for the generated fractal image, related to matrix $M$.
\end{rem}
	
Note: a matrix declared an ``initial" has the order 1, even if it is, in fact, a resultant order $N$ matrix generated previously.\\

For clarity, a simple notation is used to stress the name of the fractal (using abbreviation) and its order (in the form `oN'). 
E.g., for the Sierpinski triangle fractal (STF) order 7 -- both generated matrix and plotted picture would be denoted as STFo7.

Here is another important definition:
\begin{defn}
		A matrix containing zeros and ones will be called the \textbf{fractal matrix}
		if it has at least one zero and one number one. 
\end{defn}

Some initial matrices (even fractal ones) can't, produce fractals. So, here is a definition for them:
\begin{defn}
	A matrix is called \textbf{degraded} if it contains all zeros or only one number one.
\end{defn}

To generate and plot KPB fractals R helper functions from [3] will be use. It should be explained that R helper 
functions just simplifying and standardizing generating and plotting for fractal matrices. But, in fact, R 
language has a built-in operator $\%x\%$ for the Kronecker product, i.e., for any matrix $M$ it could be used as $M\%x\%M$.

To start, two \emph{Sierpinski canonical fractals} --  triangle and carpet -- are shown in the figures below. 

\begin{figure}[!ht]
	\begin{minipage}[b]{0.5\textwidth}
		\includegraphics[width=\textwidth]{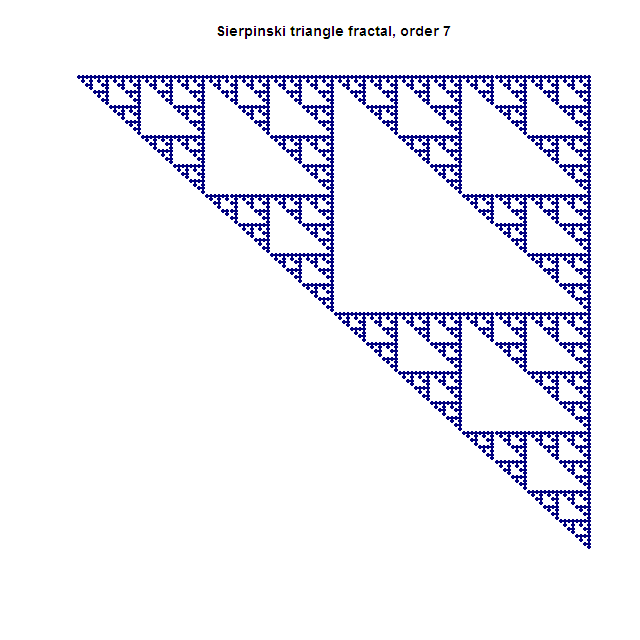}
		\caption{STF, order 7}
	\end{minipage}%
	\begin{minipage}[b]{0.5\textwidth}
		\includegraphics[width=\textwidth]{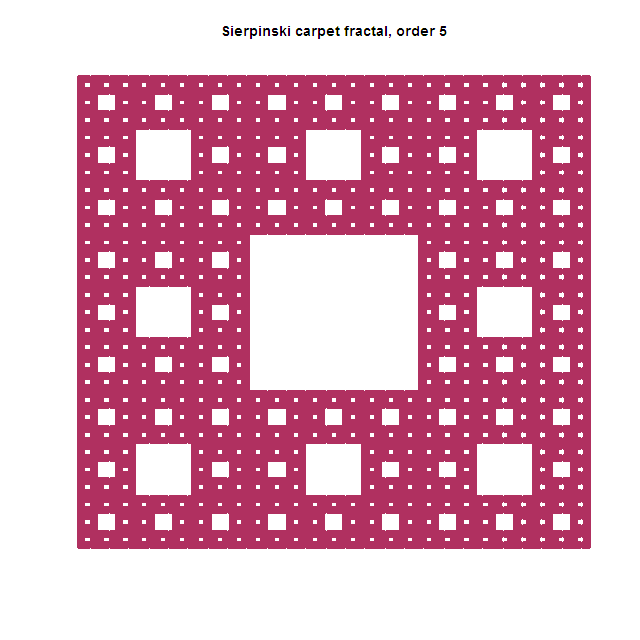}
		\caption{SCF, order 5}
	\end{minipage}%
\end{figure}

The initial fractal matrices for these two fractals and how to generate and plot them using 
helper functions from [3] is shown in the R script \#1 below.

\begin{exm} R script \#1:\\ 
\begin{verbatim}
    #### R SCRIPT BEGIN
    ## Initial matrix for Sierpinski triangle fractal (STF)
    STFo1 <- matrix(c(1,1, 0,1), nrow=2, ncol=2, byrow=TRUE);
    ## Initial matrix for Sierpinski carpet fractal (SCF)
    SCFo1 <- matrix(c(1,1,1,  1,0,1, (1,1,1), nrow=3, ncol=3, byrow=TRUE);
    ## Generate and plot fractals
    gpKronFractal(STFo1, 7, "STFo7", "navy", "Sierpinski triangle fractal");
    gpKronFractal(SCFo1, 5, "SCFo5", "maroon", "Sierpinski carpet fractal");
    #### R SCRIPT END
\end{verbatim}
\end{exm} 
	
Now, only matrices of the second order will be built (without plotting) in R.
\begin{exm} R script \#2 and its output:\\ 
\begin{verbatim}
    #### R SCRIPT & OUTPUT BEGIN
    ## Generate 2-nd order fractal matrices
    ## Initial matrix for Sierpinski triangle fractal (STF)
    STFo1 <- matrix(c(1,1, 0,1), nrow=2, ncol=2, byrow=TRUE);
    STFo1
    STFo2 = STFo1%x%STFo1; ## self Kronecker product (order 2)
    STFo2
    ## Initial matrix for Sierpinski carpet fractal (SCF)
    SCFo1 <- matrix(c(1,1,1,  1,0,1, 1,1,1), nrow=3, ncol=3, byrow=TRUE);
    SCFo1
    SCFo2 = SCFo1%x%SCFo1; ## self Kronecker product (order 2) 
    SCFo2
    
    ## OUTPUT BEGIN:
    > ## Generate 2-nd order fractal matrices
    > ## Initial matrix for Sierpinski triangle fractal (STF)
    > STFo1 <- matrix(c(1,1, 0,1), nrow=2, ncol=2, byrow=TRUE);
    > STFo1
    [,1] [,2]
    [1,]    1    1
    [2,]    0    1
    > STFo2 = STFo1%x%STFo1; ## self Kronecker product (order 2)
    > STFo2
    [,1] [,2] [,3] [,4]
    [1,]    1    1    1    1
    [2,]    0    1    0    1
    [3,]    0    0    1    1
    [4,]    0    0    0    1
    > ## Initial matrix for Sierpinski carpet fractal (SCF)
    > SCFo1 <- matrix(c(1,1,1,  1,0,1, 1,1,1), nrow=3, ncol=3, byrow=TRUE);
    > SCFo1
    [,1] [,2] [,3]
    [1,]    1    1    1
    [2,]    1    0    1
    [3,]    1    1    1
    > SCFo2 = SCFo1%x%SCFo1; ## self Kronecker product (order 2) 
    > SCFo2
    [,1] [,2] [,3] [,4] [,5] [,6] [,7] [,8] [,9]
    [1,]    1    1    1    1    1    1    1    1    1
    [2,]    1    0    1    1    0    1    1    0    1
    [3,]    1    1    1    1    1    1    1    1    1
    [4,]    1    1    1    0    0    0    1    1    1
    [5,]    1    0    1    0    0    0    1    0    1
    [6,]    1    1    1    0    0    0    1    1    1
    [7,]    1    1    1    1    1    1    1    1    1
    [8,]    1    0    1    1    0    1    1    0    1
    [9,]    1    1    1    1    1    1    1    1    1
    #### R SCRIPT & OUTPUT END
\end{verbatim}
\end{exm} 

Although R script output is very clear, in the Fig. 3, 4 below find the initial, second and third order matrices 
plotted using big square ``dots". This makes it even easier to determine the basic figure, counting and scaling.

\begin{figure}[!ht]
	\centering
	\begin{minipage}[b]{0.49\textwidth}
		\includegraphics[width=\textwidth]{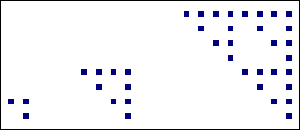}
		\caption{Schemes of order 1, 2 and 3 matrices for STF}
	\end{minipage}%
	\begin{minipage}[b]{0.49\textwidth}
	\includegraphics[width=\textwidth]{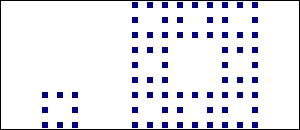}
	\caption{Schemes of order 1 and 2 matrices for SCF}
\end{minipage}%
\end{figure}


\section{Calculating Hausdorff-Besicovitch dimension}

A very simple formula to calculate Hausdorff-Besicovitch dimension of the fractal $F$ was introduced in [5] and further 
explain in [6,7] and it fits our goal:
\[
\HausB F = log N / log S,
\]
 
\noindent where $N$ - is the number of self-similar objects; $S$ - is the scale factor.

Let's apply the above formula to the Sierpinski triangle fractal.
Looking at Fig. 3 we can find that in the order 2 matrix the number of self-similar objects (i.e., objects equal to 
the basic one in the initial matrix) $N=3$,
and the scale factor $S=2$. So, $\HausB STF = log3/log2 = 1.584963$.

Applying the same formula to the Sierpinski carpet fractal (see Fig. 4) we have: the number of self-similar 
objects $N=8$, the scale factor $S=3$ and\\ $\HausB SCF = log8/log3 = 1.892789$.
  
As result of intensive testing, it was discovered that in case of KPB fractals:
\begin{itemize}
  \item $N$ is the ratio $d2/d1$, where $d2$ is a number of dots in the resultant matrix of the second order 
and $d1$ is a number of dots in the initial matrix.
  \item $S$ is the ratio $m2/m1$, where $m2$ is a number of rows in the resultant matrix of the second order
and $m1$ is a number of rows in the initial matrix.
\end{itemize}

The new simple formula to calculate Hausdorff-Besicovitch dimension of the fractal $F$ is the following: \\ 
\[
\HausB F = log(d2/d1) / log(m2/m1) ,
\]

\noindent where $d1, d2, m1, m2$ - are,  accordingly, numbers of dots and rows in the initial and the second order matrices.

Such amazing simplicity could be explained by the fact that the Kronecker product is building resultant block-matrix in a uniform manner, so,
these two ratios are always describing correctly both the number of self-similar objects and the scale factor.
Moreover, these ratios are the same between other orders, e.g., between the second and third orders, etc.

Based on this discovered formula, -- a very simple R script was created. Just 7 lines of R code:

\emph{R script dimHB4kpf.R}
\begin{verbatim}
    ## dimHB4kpf.R.txt 3/11/17 aev
    ## Hausdorff-Besicovitch dimension for Kronecker product based
    ## fractals
    ## Note: only for Kronecker product based fractals created from
    ## the fractal matrix!
    ## dimHB4kpf(mat, ttl) -
    ## where: mat - initial matrix (filled with 0/1); ttl - title.
    dimHB4kpf <- function(mat, ttl="") {
      m1 = nrow(mat); dn1 = sum(mat!=0);
      matr = mat%x%mat; ##self Kronecker product 
      m2 = nrow(matr); dn2 = sum(matr!=0);
      dimHB = log(dn2/dn1, m2/m1);
      cat(" *** dimHB:", dimHB, ttl, "\n");
    }
\end{verbatim}

Presented below testing R script \#3 is pretty simple. It includes two lines of code for
each of 7 selected well-known fractals.

\begin{exm} R script \#3:\\
	\begin{verbatim}
    ## Testing dimHB for 7 well-known fractals:
    STF <- matrix(c(1,1, 0,1), ncol=2, nrow=2, byrow=TRUE);
    dimHB4kpf(STF, "'Sierpinski triangle'")
    SCF <- matrix(c(1,1,1, 1,0,1, 1,1,1), ncol=3, nrow=3, byrow=TRUE);
    dimHB4kpf(SCF, "'Sierpinski carpet'")
    PTFm3 <- matrix(c(1,1,1, 0,1,1, 0,0,1), ncol=3, nrow=3, byrow=TRUE);
    dimHB4kpf(PTFm3, "'Pascal triangle modulo 3'")
    PTFm5 <- matrix(c(1,1,1,1,1, 0,1,1,1,1, 0,0,1,1,1, 0,0,0,1,1, 0,0,0,0,1), 
    ncol=5, nrow=5, byrow=TRUE);
    dimHB4kpf(PTFm5, "'Pascal triangle modulo 5'")
    VF <- matrix(c(0,1,0, 1,1,1, 0,1,0), ncol=3, nrow=3, byrow=TRUE);
    dimHB4kpf(VF, "'Vicsek fractal'")
    HGF <- matrix(c(1,1,0, 1,1,1, 0,1,1), ncol=3, nrow=3, byrow=TRUE);
    dimHB4kpf(HGF, "'Hexagon/Hexaflake fractal'")
    BF <- matrix(c(1,0,1, 0,1,0, 1,0,1), ncol=3, nrow=3, byrow=TRUE);
    dimHB4kpf(BF, "'Box fractal'")
	\end{verbatim}
\end{exm}

Note: last 5 fractals are shown in the Fig. 5-9 below.
\newpage
\begin{figure}[!ht]
	\begin{minipage}[b]{0.5\textwidth}
		\includegraphics[width=\textwidth]{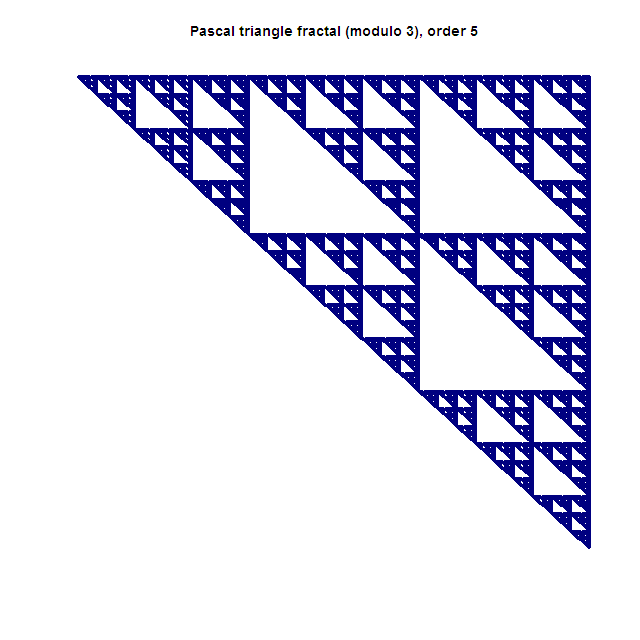}
		\caption{PTFm3o5}
	\end{minipage}%
	\begin{minipage}[b]{0.5\textwidth}
		\includegraphics[width=\textwidth]{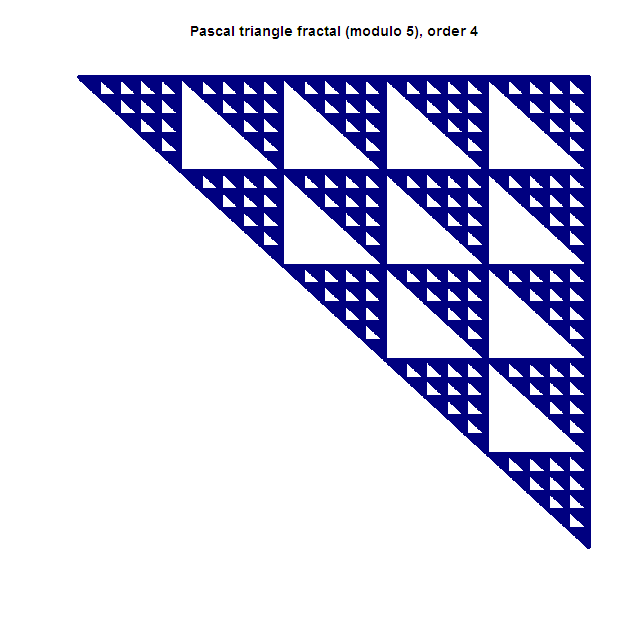}
		\caption{PTFm5o4}
	\end{minipage}%
\end{figure}

\begin{figure}[!ht]
	\begin{minipage}[b]{0.33\textwidth}
		\includegraphics[width=\textwidth]{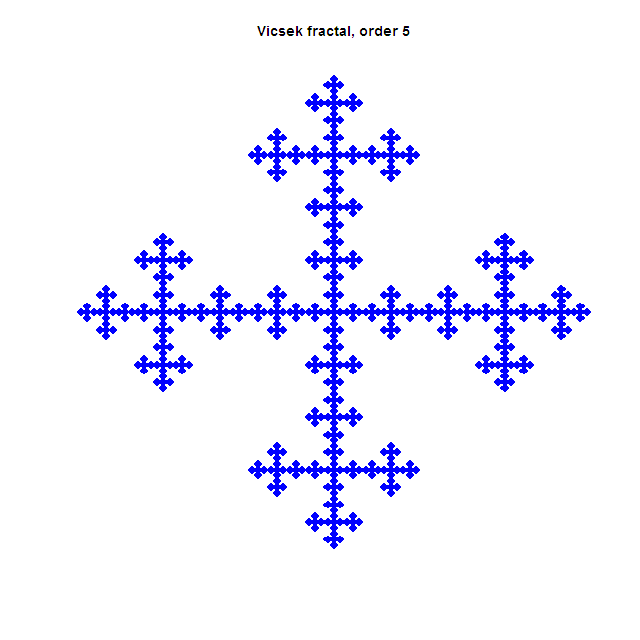}
		\caption{VFo5}
	\end{minipage}%
	\begin{minipage}[b]{0.33\textwidth}
		\includegraphics[width=\textwidth]{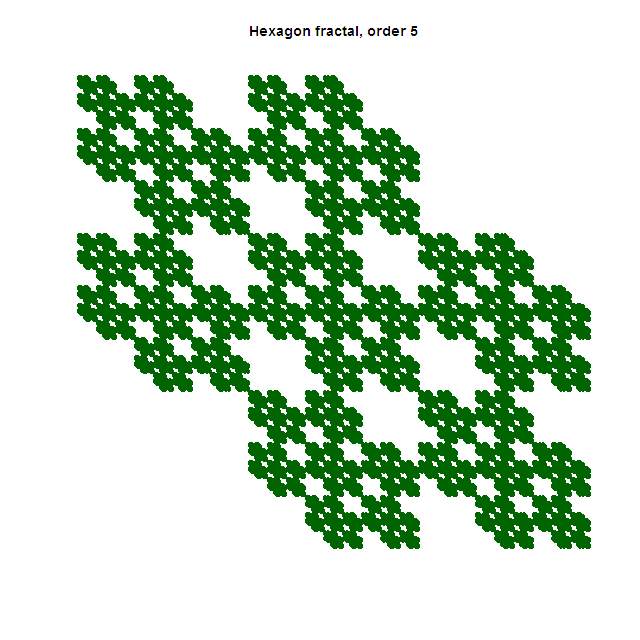}
		\caption{HGFo5}
	\end{minipage}%
	\begin{minipage}[b]{0.33\textwidth}
	\includegraphics[width=\textwidth]{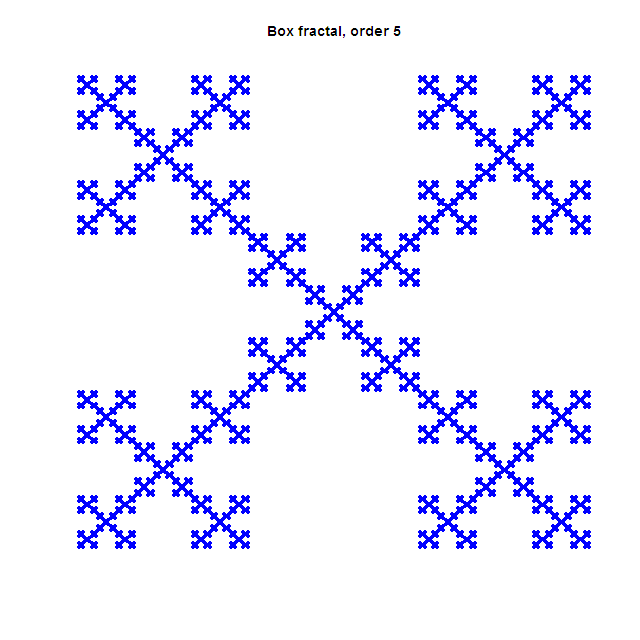}
	\caption{BFo5}
    \end{minipage}%
\end{figure}

Results of the testing script are presented below, and they are equal to results shown in [8] for 
the same fractals.

\begin{exm} R script \#3 output:\\
\begin{verbatim}
    ## OUTPUT BEGIN:
    > ## Testing dimHB for 7 well-known fractals:
    > STF <- matrix(c(1,1, 0,1), ncol=2, nrow=2, byrow=TRUE);
    > dimHB4kpf(STF, "'Sierpinski triangle'")
    *** dimHB: 1.584963 'Sierpinski triangle' 
    > SCF <- matrix(c(1,1,1, 1,0,1, 1,1,1), ncol=3, nrow=3, byrow=TRUE);
    > dimHB4kpf(SCF, "'Sierpinski carpet'")
    *** dimHB: 1.892789 'Sierpinski carpet' 
    > PTFm3 <- matrix(c(1,1,1, 0,1,1, 0,0,1), ncol=3, nrow=3, byrow=TRUE);
    > dimHB4kpf(PTFm3, "'Pascal triangle modulo 3'")
    *** dimHB: 1.63093 'Pascal triangle modulo 3' 
    > PTFm5 <- matrix(c(1,1,1,1,1, 0,1,1,1,1, 0,0,1,1,1, 0,0,0,1,1, 
    +        0,0,0,0,1), ncol=5, nrow=5, byrow=TRUE);
    > dimHB4kpf(PTFm5, "'Pascal triangle modulo 5'")
    *** dimHB: 1.682606 'Pascal triangle modulo 5' 
    > VF <- matrix(c(0,1,0, 1,1,1, 0,1,0), ncol=3, nrow=3, byrow=TRUE);
    > dimHB4kpf(VF, "'Vicsek fractal'")
    *** dimHB: 1.464974 'Vicsek fractal' 
    > HGF <- matrix(c(1,1,0, 1,1,1, 0,1,1), ncol=3, nrow=3, byrow=TRUE);
    > dimHB4kpf(HGF, "'Hexagon/Hexaflake fractal'")
    *** dimHB: 1.771244 'Hexagon/Hexaflake fractal' 
    > BF <- matrix(c(1,0,1, 0,1,0, 1,0,1), ncol=3, nrow=3, byrow=TRUE);
    > dimHB4kpf(BF, "'Box fractal'")
    *** dimHB: 1.464974 'Box fractal' 
    ## OUTPUT END
\end{verbatim}
\end{exm}

\section{Conclusion}

The discovered new simple formula to calculate the Hausdorff-Besicovitch dimension together with 
offered compact R script (dimHB4kpf.R) can help studies related to KPB fractals and even to any fractals in general.

For example, if any fractal can be simulated using matrix presentation and the Kronecker power, 
then its Hausdorff-Besicovitch dimension can be calculated instantly, - helping understand many peculiarities of 
this fractal.

Actually, it was done already here for Hexagon fractal (see Fig. 8). Our plotted Hexagon fractal is very similar to 
the Hexaflake in [8], but it looks rather as a distorted Hexaflake. Anyway, they both have the equal dimension.

Here is another appropriate and interesting example related to the Rauzy gasket (RG) described in [9]. The authors of [9] 
proved ``that the Rauzy gasket is homeomorphic to the usual Sierpinski gasket". If compare the Fig. 1 in the mentioned 
article (which is presented the Rauzy gasket) and the Fig. 1 above (presented Sierpinski triangle fractal, aka 
Sierpinski gasket), - it is very clear that the RG is looking like a distorted version of STF, and it is actually copying 
the structure of STF.

It was already calculated that $\HausB STF = 1.584963$, so  $\HausB RG$  
is the same (at least, expected to be almost the same).
This value is in accordance with the fact proved in [10]: ``the Hausdorff dimension of the Rauzy gasket is less than 2".
One can look at this value as a reasonable approximation, until the precise value would be found, or it would be proved 
that dimensions of both gaskets are equal. 



\end{document}